\documentclass[10pt,twoside]{classe-ma_e}
\usepackage{amssymb,amsbsy,amsmath,amsfonts,amssymb,amscd}
\usepackage{latexsym,euscript,exscale,epic,eepic,epsfig}
\usepackage[english,francais]{babel}
\usepackage{times}
\DeclareFontFamily{U}{rsf}{}
\DeclareFontShape{U}{rsf}{m}{n}{
  <5> <6> rsfs5 <7> <8> <9> rsfs7 <10-> rsfs10}{}
\DeclareMathAlphabet{\mathscr}{U}{rsf}{m}{n}

\DeclareMathAlphabet{\mathgth}{U}{euf}{m}{n}

\DeclareFontFamily{U}{cyr}{}
\DeclareFontShape{U}{cyr}{m}{n}{
  <5> wncyr5 <6> wncyr6 <7> wncyr7 <8> wncyr8 <9> wncyr9 <10-> wncyr10}{}
\DeclareMathAlphabet{\mathcyr}{U}{cyr}{m}{n}

\input cyracc.def

\newcommand{\cE}{{\mathscr E}}
\newcommand{\cF}{{\mathscr F}}
\newcommand{\cG}{{\mathscr G}}

\newcommand{\cL}{{\mathscr L}}

\newcommand{\cO}{{\mathscr O}}

\newcommand{\cU}{{\mathscr U}}

\newcommand{\FMXM}{\Phi_{X\ra M}}
\newcommand{\FMMX}{\Phi_{M\ra X}}

\newcommand{\D}{{\mathbf D}_{\mathrm{coh}}^b}

\newcommand{\chk}{{\scriptscriptstyle\vee}}
\newcommand{\R}{\mathbf{R}}
\newcommand{\Ld}{\mathbf{L}}
\newcommand{\lotimes}{\stackrel{\Ld}{\otimes}}

\DeclareMathOperator{\Jac}{Jac}

\DeclareMathOperator{\Hom}{Hom}

\DeclareMathOperator{\Pic}{Pic}

\DeclareMathOperator{\rk}{rk}
\DeclareMathOperator{\Ext}{Ext}

\newcommand{\ra}{\rightarrow}

\newcommand{\C}{\mathbf{C}}

\newcommand{\Z}{\mathbf{Z}}

\newcommand{\iso}{\cong}

\newcommand{\pj}{\mathbf{P}}

\renewcommand{\phi}{\varphi}

\newtheorem{theorem}{Theorem}

\newtheorem{e-proposition}[theorem]{Proposition}
\newtheorem{corollary}[theorem]{Corollary}
\newtheorem{e-definition}[theorem]{Definition\rm}

\ComParit{S1631-073X}
\PIT{FLA}
\PXMA{02290}
\Add{2}
\Volume{334}
\Year{2002}
\FirstPage{469}
\LastPage{472}
\AuteurCourant{Andrei C\u ALD\u ARARU}
\TitreCourant{Fiberwise stable bundles on elliptic threefolds} 
\Journal
\Rubrique{G\'eom\'etrie alg\'ebrique}{Algebraic Geometry}
\SousRubrique{}{}
\PresentePar{Jean-Pierre}{Demailly}
\Recu{29 November 2001}{}{7 January 2002}
\begin{document}
\selectlanguage{english}
\title{%
Fiberwise stable bundles on elliptic threefolds with relative Picard number one
}
\author{%
Andrei C\u ALD\u ARARU
}
\address{%
Mathematics Department, University of Massachusetts, Amherst, MA
01003, USA \\ 
E-mail: andreic@math.umass.edu
}
\maketitle

\thispagestyle{empty}
\begin{Abstract}{%
We show that fiberwise stable vector bundles are preserved by relative
Fourier-Mukai transforms between elliptic threefolds with relative
Picard number one.  Using these bundles we define new invariants of
elliptic fibrations, and we relate the invariants of a space with
those of a relative moduli space of stable sheaves on it.  As a
byproduct, we calculate the intersection form of a certain new example
of an elliptic Calabi-Yau threefold.  }\end{Abstract}
\selectlanguage{french}
\begin{Ftitle}{%
Fibr\'es vectoriels relativement stables sur vari\'et\'es elliptiques
de dimension trois dont le num\'ero relatif de Picard est un
}\end{Ftitle}
\begin{Resume}{%
Nous prouvons que les fibr\'es vectoriels relativement stables sont
pr\'eserv\'es par des transform\'ees de Fourier-Mukai entre vari\'et\'es
elliptiques de dimension trois dont le num\'ero relatif de Picard est
un.  En utilisant ces fibr\'es nous d\'efinissons des nouveaux
invariants de vari\'et\'es elliptiques, et nous \'etudions la
relation entre les invariants d'une vari\'et\'e et ceux d'un \'espace
relatif de modules des fibr\'es stables sur elle.  Ces r\'esultats
nous permettent de calculer la forme d'intersection sur un certain
nouvel exemple de vari\'et\'e de Calabi-Yau de dimension
trois. }\end{Resume}

%
\setcounter{section}{0}
\selectlanguage{english}

\section*{Introduction}
The object of this note is to generalize to the case of elliptic
threefolds with relative Picard number one classic results regarding
stable vector bundles on elliptic curves and their Fourier-Mukai
transforms.  Several results of this type are known: partial results
by Bridgeland~(\cite{BriEll}) for arbitrary elliptic fibrations, and
strong results by Bartocci, Bruzzo, Hern\'andez-Ruip\'erez and
Mu\~noz-Porras~(\cite{BBHRP}) for elliptic fibrations without
reducible fibers.  Our primary interest is in applying this theory to
the study of elliptic Calabi-Yau threefolds, so we need to have
results general enough to handle reducible curves of arithmetic genus
1.

Stable bundles on an elliptic curve $E$, of rank $r$ and degree $d$,
$(r, d) = 1$, are well understood:
\smallskip

\textsc{Theorem} (Atiyah~\cite{Ati}).\;\;--\;\; 
{\it For every line bundle $\cL$ of degree $d$ on $E$, there exists a
unique stable vector bundle $V(r, \cL)$ of rank $r$ and determinant
$\cL$.}
\medskip

\textsc{Theorem} (Tu~\cite{Tu}).\;\;--\;\; 
{\it A vector bundle $\cF$ on $E$ of rank $r$ and degree $d$ is stable
if and only if it is simple (i.e.\ $\Hom(\cF, \cF) = \C$).}
\medskip

Let $a,b$ be coprime integers, with $a>0$, and let $M$ be the moduli
space of stable bundles of rank $a$, degree $b$ on $E$ (which is a
fine moduli space).  Fix a universal sheaf $\cU$ on $X\times M$, and
let $c$ be the degree of the restriction $\cU|_{\{x\}\times M}$ for
$x\in X$, which is independent of $x$.  Let $p_E,\ p_M$ be the
projections from $E\times M$ to $E$, $M$, respectively.
\medskip

\textsc{Theorem} (Bridgeland~\cite{BriEll}).\;\;--\;\; 
{\it Let $\cF$ be a stable vector bundle on $M$, of slope $\mu(\cF)
\neq -c/a$.  Then the Fourier-Mukai transform $\cG$ of $\cF$, $\cG =
\Phi_{M\ra E}^\cU(\cF) = \R p_{E,*}(p_M^*(\cF) \otimes \cU)$ is a
stable vector bundle on $E$ of slope $\mu(\cG) \neq b/a$, possibly
shifted.}
\medskip

The stability of $\cF$ is used in two ways in the proof of the above
theorem.  It is first used to conclude that $\cG$ consists of a single
vector bundle, possibly shifted.  Then, since $\cF$ is simple and
simplicity is preserved by Fourier-Mukai transforms, $\cG$ is simple,
and thus stable by Tu's theorem.

\section*{Fiberwise stable bundles}
Our goal is to replace the elliptic curve $E$ in the above results by
an elliptic threefold $X/S$, with $X$ and $S$ smooth, thus allowing
for some reducible fibers.  In such generality, only partial results
are known, mainly due to Bridgeland~(\cite{BriK3}).  Note that the
theorems of Atiyah and Tu hold when $E$ is replaced by an irreducible
curve of arithmetic genus 1, but fail if $E$ is reducible.  We'll
assume that the fibration $X/S$ has no multiple fibers, its relative
Picard number $\rho(X/S) = \rho(X)-\rho(S)$ is equal to one, $\Pic(X)$
has no torsion, and the general singular fiber of $X/S$ is
irreducible.  These conditions are often satisfied by general members
of families of elliptic Calabi-Yau threefolds~(\cite{CalEll}).  The
assumption that $X$ and $S$ are smooth ensures that all the fibers of
$X\ra S$ are Cohen-Macaulay.

Given such an elliptic fibration $X/S$, and $r$, $d$ coprime integers
with $r>0$, denote by $V_X(r,d)$ the class of vector bundles of rank
$r$ on $X$, whose restriction to each fiber of $X/S$ is stable of
degree $d$.  Since we are dealing with possibly reducible fibers, we
use Gieseker's definition of stability.  In general, the definition of
$V_X(r,d)$ depends on the choice of polarization in each fiber, but
the assumption $\rho(X/S) = 1$ implies that each fiber is polarized in
a unique way (up to multiples) by the restriction of a polarization of
$X$.  From here on, we'll always assume the fibers polarized in this
way.  The case $r<0$ can also be included in the definition, by
considering elements of $V_X(r,d)$ as objects in $\D(X)$ which consist
of a single sheaf, possibly shifted (if $\cE$ is a $V(r,d)$, then
$\cE[1]$ is a $V(-r, -d)$).

\begin{e-proposition}
\label{prop:1}
Elements of $V_X(r,d)$ differ by tensoring by pull-backs of line
bundles in $\Pic(S)$ and by even shifts in the derived category.
\end{e-proposition}

\begin{proof}
Follows from Atiyah's result on irreducible curves of genus 1 and
Hartogs' theorem.
\end{proof}

We want to study the behaviour of the $V(r,d)$-bundles under
Fourier-Mukai transforms.  Let $n$ be the smallest positive degree of
a multisection of $X/S$ (alternatively, this is the smallest degree of
a polarization of the fibers which is the restriction of a
polarization from $X/S$).  Consider integers $a>0$ and $b$ so that
$(na,b) = 1$, and let $M/S$ be the relative moduli space of stable
sheaves of rank $a$, degree $b$ on the fibers of $X/S$, in the sense
of Simpson.  The space $M$ has a natural map to $S$ which makes it
into an elliptic fibration which satisfies all the conditions imposed
on $X/S$.  The fibration $M/S$ is a fine moduli space, and the
extension by zero $\cU$ to $X\times M$ of a universal sheaf on
$X\times_S M$ induces a Fourier-Mukai transform $\FMMX^{\cU}:\D(M) \ra
\D(X)$~(\cite{BriEll},~\cite{BriK3}).  Let $c$ be the degree of the
restriction of $\cU$ to a fiber $\{x\} \times M$ for $x\in X$ and
$e=(bc-1)/a$.  For $r$ and $d$ coprime define $V_M(r,d)$ in an
analogous fashion to $V_X(r,d)$.

\begin{theorem}
\label{thm:1}
If $\cF$ is a $V_M(r,d)$, for $r,d$ coprime, $d/r\neq -c/a$, then the
Fourier-Mukai transform $\cG = \FMMX^{\cU}(\cF)$ is a $V_X(r',d')$,
where $(r', d')$ are given by 
\[ \left (\!\!\!\begin{array}{l}r' \\ d'\end{array}\!\!\!\right) = 
\left (\!\!\!\begin{array}{ll}c & a \\ e & b\end{array}\!\!\!\right )
\left (\!\!\!\begin{array}{l}r \\ d\end{array}\!\!\!\right). \] 
\end{theorem}

The relationship between $(r', d')$ and $(r,d)$ has been known
previously~(\cite{BriEll}); what is new is the fact that $\cG$ is a
$V_X(r', d')$.  The proof of Theorem~\ref{thm:1} follows the same
lines as the proof of Bridgeland's theorem, and in order to conclude
that $\cG$ is a single, locally free sheaf on $X$ we need to have the
restriction of $\cF$ to every fiber of $M/S$ be stable.  (This part of
the proof goes through with almost no restrictions on $X/S$.)  But
although we are guaranteed that the restriction of $\cG$ to all the
fibers of $X/S$ is simple, we can no longer apply Tu's result to
conclude that it is stable.  Thus, for general fibrations, the
situation is asymmetric (we need to start with a sheaf that is {\em
stable} on all the fibers, and we end up with one that is only {\em
simple}).  The crucial observation in the case of relative Picard
number one is the following:

\begin{e-proposition}
\label{prop:2}
Let $E$ be a (possibly reducible) Cohen-Macaulay curve of arithmetic
genus 1, and let $\cE$ be a locally free sheaf on $E$ whose
determinant is either ample or anti-ample.  Then $\cE$ is simple if
and only if it $\cE$ is stable with respect to the polarization
induced by $\det\cE$ or $-\det\cE$.
\end{e-proposition}

\begin{proof}
One implication is trivial, so assume that $\cE$ is simple, and that
$\cL = \det\cE$ is ample.  Then $\cL = \cO(\sum a_i P_i)$, with
$a_i>0$ and $P_i\in E^{\mathrm{smooth}}$.  For any torsion-free sheaf
$\cF$ on $E$, let
\[ \mu_\cL(\cF) = \frac{\chi(\cF)}{\sum a_i \rk_{P_i}(\cF)}, \]
where $\rk_{P_i}(\cF) = \dim \cF_{P_i}$.  The reduced Hilbert
polynomial of $\cF$, computed with respect to $\cL$, is equal to
$t+\mu_\cL(\cF)$.  Thus, to show that $\cE$ is stable with respect to
$\cL$, we need to show that if $\cF$ is a proper subsheaf of $\cE$,
then $\mu_\cL(\cF) < \mu_\cL(\cE)$.  A straightforward computation
shows that $\mu_\cL(\cF\otimes \cE^\chk) = \mu_\cL(\cF) - \mu_\cL(\cE)$.
If $\cF$ is a subsheaf of $\cE$, $\Hom(\cF, \cE) \neq 0$ and thus
$H^1(E, \cF\otimes \cE^\chk) = \Ext^1(\cE, \cF) = \Hom(\cF, \cE) \neq
0$ by Serre duality (which can be applied because $\cE$ is locally
free).  If $\mu_\cL(\cF) \geq \mu_\cL(\cE)$, then $\chi(\cF\otimes
\cE^\chk) \geq 0$, so we conclude that $\Hom(\cE, \cF) = H^0(E,
\cF\otimes\cE^\chk) \neq 0$.  Thus there is a non-zero map $\cE\ra
\cF$, which composed with the inclusion $\cF\subset \cE$ yields a
non-trivial map $\cE\ra\cE$, contradicting the assumption that $\cE$
is simple.  The case $\det\cE$ anti-ample is treated in a similar way.
\end{proof}

\begin{corollary}
Let $\cE$ be a locally free sheaf on $X$.  Then its restriction to a
fiber of $X/S$ is stable (with respect to the unique polarization of
the fiber coming from $X$) if and only if it is simple.
\end{corollary}

Note that $\det\cE$ must restrict to a multiple of the polarization
$\cO_{X/S}(1)$ on each fiber, because $\rho(X/S) = 1$.  If this
restriction is non-trivial, we can apply Proposition~\ref{prop:1}.  If
$\det\cE$ is the trivial line bundle in each fiber, $\cE(1)$ is stable in
each fiber with respect to $\cO_{X/S}(1)$, which implies that $\cE$
is fiberwise stable.

\section*{Invariants of elliptic fibrations}
We'd like to have invariants that enable us to compare $X$ and $M$.
Define $P_X(r,d)$ to be the class of elements of $\D(S)$ that are of
the form $\R\pi_{X,*}\cF$ for $\cF$ a $V_X(r,d)$ ($\pi_X:X\ra S$ is
the structural map of the fibration $X/S$).  Define $P_M(r,d)$
analogously.  

The following two propositions are applications of
Proposition~\ref{prop:1}, the projection formula, Grothendieck-Serre
duality and the fact that the dual of a $V(r,d)$ is a $V(r,-d)$, which
follows from Proposition~\ref{prop:2}.

\begin{e-proposition}
\label{prop:3}
For $d\neq 0$, elements of $P_X(r,d)$ are vector bundles (possibly
shifted) which differ by tensoring by line bundles in $\Pic(S)$ and by
even shifts in the derived category.
\end{e-proposition}

\begin{e-proposition}
\label{prop:5}
$P_X(r,d) = P_X(r, -d)^\chk[1] = P_X(-r, d)^\chk = P_X(-r, -d)[1]$,
where $P_X(r, -d)^\chk[1]$ is obtained by dualizing and shifting every
element of $P_X(r,-d)$, and $P_X(-r, d)^\chk$, $P_X(-r, -d)[1]$ are
defined in a similar way.
\end{e-proposition}

\begin{theorem}
\label{thm:3}
We have $P_X(r',d') = \R\pi_{M,*} (V_M(r,d) \otimes V_M(b,e))$, where
$a, b, c, e$ and $r, d, r', d'$ are as in Theorem~\ref{thm:1}.  In
particular, if $r=1$, $P_X(c+ad, e+bd) = P_M(b, e+bd)$.
\end{theorem}

{\it Proof\/.~~--~\kern.3em}\ignorespaces
Let $\cF$ be a $V_M(r,d)$, and let $\cG =
\FMMX^\cU(\cF) = \R p_{X,*}(p_M^* \cF \lotimes \cU)$, where $p_X, p_M$
are the projections from $X\times M$ to $X$ and $M$, respectively.
Then $\cG$ is a $V_X(r', d')$, $\R\pi_{X,*}\cG = \R\pi_{M,*}(\cF
\otimes\R p_{M,*}\cU)$, and all we need to show is that $\R
p_{M,*}\cU$ is a $V_M(b,e)$.  But $\R p_{M,*} \cU = \FMXM^\cU(\cO_X)$,
so the result follows from Theorem~\ref{thm:1} once we show that $X$
is the relative moduli space of stable sheaves of rank $a$, degree $c$
on $M$, and $\cU$ is a universal sheaf for this moduli problem.  The
proof of this fact (which is known for elliptic
surfaces,~\cite{BriEll}) will be presented elsewhere.
\qed\goodbreak\vskip6pt

\begin{theorem}
\label{thm:2}
If $X/S$ also satisfies the property of being generic
(see~\cite{CalEll} for a definition), then $P_X(r,d) = P_X(r',d)$ for
$r = r'\bmod d$.
\end{theorem}
\smallskip

The proof of Theorem~\ref{thm:2} is more involved, requiring the use
of twisted derived categories~(\cite{Cal}), and will be included in a
future paper.  In the case of a single elliptic curve $E$, the proof
does not rely on twisted sheaves, and we include a sketch here.  If
$\cF$ is a $V_E(r,d)$, then its transform $\cG$ to $\Jac(E)$ is a
$V_{\Jac(E)}(-d, r)$.  On $\Jac(E)$ there is a naturally defined
$\cO_{\Jac(E)}(-1)$, and $\cG(-1)$ is a $V_{\Jac(E)}(-d, r+d)$.
Transforming $\cG(-1)$ back to $E$ yields a $V_E(r+d, d)$, whose
global sections can be shown to be the same as those of $\cF$.

\section*{An explicit calculation}
In~\cite[6.2.2]{Cal} we studied an explicit generic elliptic
Calabi-Yau threefold $X/\pj^2$, constructed as the pfaffian of a
certain $5\times 5$ matrix of bihomogeneous forms on
$\pj^2\times\pj^4$. We described all the minimal birational models of
$X$, which consist of $X$, a flop $X^+$ of $X$ which is contained in
$\pj^4\times\pj^5$ and has no elliptic fibration structure, and a flop
$X^{++}$ of $X^+$ contained in $\pj^5\times\pj^2$, which has an
elliptic fibration structure given by the map $X^{++}\ra \pj^2$.
There is no apparent relation between the original elliptic fibration
structure on $X$ and the one on $X^{++}$.  Since $X^{++}$ is obtained
from $X$ through a sequence of flops, $\D(X^{++}) \iso
\D(X)$~(\cite{BriFlops}).

On $X$ we have $n = 5$; taking $a=1$, $b=2$ we obtain an elliptic
Calabi-Yau threefold $M/\pj^2$, with $\D(X)\iso\D(M)$.  It can be
shown that $V_X(1,5)$ and $V_X(3,5)$ are non-empty, and thus $V_X(r,
5)$ is non-empty for all $r$ relatively prime to 5.  Comparing the
invariants $P_X(r, 5)$ and $P_M(r,5)$ for all $r$ relatively prime to
5 and using Grothendieck-Riemann-Roch, we compute $c_2(M)$ and the
cubic form on $H^2(M, \Z)$.  These topological invariants can be
computed for $X^{++}$, and they are the same as those of $M$.

Thus the two Calabi-Yau threefolds $X^{++}$ and $M$ have $\D(X^{++})
\iso \D(M)$, and have the same $c_2$ and cubic form.  In view of
Wall's results~(\cite{Wal}), it seems natural to conjecture that
$X^{++}$ and $M$ are deformation equivalent.  If they are
non-isomorphic, this is likely to lead to a new counterexample to
Torelli for Calabi-Yau threefolds (the equivalence between the derived
categories of $X$ and $M$ induces an isomorphism between their
polarized Hodge structures).  In this counterexample, $X^{++}$ and $M$
would be non-birational, since we can enumerate all the birational
models of $X^{++}$ and none of them is isomorphic to $M$.  If $X^{++}$
and $M$ are isomorphic, this would suggest that there is some deeper
phenomenon behind this occurrence, which is worth investigating.  It
would also yield a new automorphism of the derived category of
$X^{++}$, which is interesting to study in view of implications to Kontsevich's
homological mirror symmetry conjecture.

\Acknowledgements{I am grateful to Mark Gross for suggesting me to
look at elliptic threefolds for counterexamples to Torelli, and to
Sorin Popescu, Titus Teodorescu and Eyal Markman for numerous useful
discussions.}

\end{document}